\documentclass[12pt]{article}
\usepackage{latexsym, amsmath, amssymb, amsfonts, amscd}
\usepackage{color}
\usepackage{graphics}
\usepackage[dvips]{epsfig}

\newtheorem{theorem}{Theorem}
\newtheorem{proposition}{Proposition}
\newtheorem{lemma}{Lemma}

\newtheorem{definition}{Definition}

\newtheorem{question}{Question}

\newcommand{\calo}{{\cal O}}

\newcommand{\adL}{\mbox{\rm ad}_{\Lambda}}

\newcommand{\ad}{\mbox{\rm ad}}

\newcommand{\image}{\mbox{\rm image\ }}

\def\lcf{\lbrack\! \lbrack}
\def\rcf{\rbrack\! \rbrack}

\def\dbar{\overline\partial}

\newcommand{\CC}{\mathbb{C}}

\newcommand{\RR}{\mathbb{R}}

\def\dbar{\overline\partial}

\def\oomega{\overline\omega}
\def\Oomega{\overline\Omega}

\newcommand{\lie}[1]{\mathfrak{#1}}
\newcommand{\lieg}{\mathfrak{g}}
\newcommand{\liet}{\mathfrak{t}}
\newcommand{\liec}{\mathfrak{c}}
\newcommand{\lieh}{\mathfrak{h}}

\newcommand{\lbra}[2]{\lcf #1, #2 \rcf}

\newcommand{\classl}[1]{ [ #1 ]_{\Lambda}}
\newcommand{\classd}[1]{ [ #1 ]_{\dbar}}

\newcommand{\bproof}{\noindent{\it Proof: }}
\newcommand{\eproof}{\hfill \qed \vspace{0.2in}}
\def\qed{\rule{2.3mm}{2.3mm}}

\begin{document}
\title{\bf A Hodge-Type Decomposition of Holomorphic Poisson Cohomology
on Nilmanifolds}
\author{
Yat Sun  Poon
    and
    John Simanyi
    }
\maketitle
\begin{abstract} A cohomology theory associated to a holomorphic
Poisson structure is the hypercohomology of a bi-complex where one of the two operators is the classical
$\dbar$-operator, while the other operator is the adjoint action of the Poisson bivector with respect to
the Schouten-Nijenhuis bracket. The first page of the associated
 spectral sequence is the Dolbeault cohomology
with coefficients in the sheaf of germs of holomorphic polyvector fields.
 In this note, the authors investigate the
conditions for which this spectral sequence degenerates on the first page when the
underlying complex manifolds are nilmanifolds with an abelian complex structure.
For a particular class of holomorphic Poisson structures, this result leads to a Hodge-type
decomposition
  of the holomorphic Poisson cohomology. We provide examples when the nilmanifolds are 2-step.
\end{abstract}


\



\section{Introduction}

It is well known that complex structures and symplectic structures are examples of
generalized complex structures in the sense of Hitchin
{\cite{Marco, Hitchin-Generalized CY, Hitchin-Instanton}}.
It is now also known that a
holomorphic Poisson structure plays a fundamental role in generalized geometry
\cite{Bailey}. Since a key feature of generalized geometry is to place both complex
structures and symplectic structures within a single framework,
its deformation theory is of great interest
\cite{Goto1, Goto2, Hitchin-holomorphic Poisson}.
An understanding of
 cohomology theory of generalized geometry in general, and
 holomorphic Poisson structures in particular, becomes necessary.
On algebraic surfaces, there has been work done by Xu and collaborators
\cite{Hong-Xu, Xu}. Recently, Hong studied holomorphic Poisson cohomology on toric Poisson manifolds
\cite{Hong}.
For nilmanifolds with abelian
complex structures, there is recent work done by this author and his collaborators
\cite{CFP, CGP, GPR}.
A common feature of these works is to recognize the cohomology as the hypercohomology of a bi-complex.

In computation of hypercohomology of bi-complex,  theoretically one could apply
one of the two naturally defined spectral sequences to
complete the task. In the case of a holomorphic Poisson structure,
one of the operators is the classical $\dbar$-operator. Another operator
is the action of a holomorphic bi-vector field $\Lambda$ on the space of
polyvector fields and differential forms via the Schouten bracket.
One may choose a filtration so that the first page of
the spectral sequence consists of the Dolbeault
cohomology of a complex manifold with coefficients in the
sheaf of germs of holomorphic polyvector fields, $H^q(M, \Theta^p)$. In this case, the
$d_1$-map on the first page is the holomorphic Poisson structure via the Schouten bracket.
Another approach is to apply the holomorphic Poisson structure action first. In
such a case, the first page of the spectral sequence consists of
a holomorphic version of Lichnerowicz-Poisson cohomology \cite{Lich, Vaisman}.

We focus on the first approach, as Dolbeault cohomology is a well known classical object, and
it involves the complex structure only. It is natural to determine when this spectral
sequence degenerates quickly. In \cite{CFP, CGP},
we have seen situations where the spectral sequence
degenerates on the second page. We push our analysis from past work and investigate the possibility when
degeneracy on the first page occurs, and establish a direct relation between the holomorphic Poisson
cohomology $H^n_{\Lambda}$ and the elements in the first page of the spectral sequence.
In particular, we could formulate a result as below.

\begin{theorem}\label{corollary theorem}
Let $M$ be a nilmanifold with an abelian complex structure. Let $H^q(M, \Theta^p)$ be the
$q^{\mathit{th}}$ Dolbeault cohomology of the complex manifold with coefficients in the
sheaf of germs of holomorphic polyvector fields of degree $p$. For any invariant holomorphic Poisson
structure $\Lambda$ on $M$, let $H^n_{\Lambda}$ be the $n^{\mathit{th}}$ cohomology of the holomorphic Poisson structure. Then there exists a natural one-to-one linear map
\[
\phi: H^n_{\Lambda}\to  \oplus_{p+q=n}  H^q(M, \Theta^p).
\]
In particular,
$\dim H^n_{\Lambda}\leq \sum_{p+q=n} \dim H^q(M, \Theta^p)$.
\end{theorem}

To prove this result, we first review some basic facts on generalized complex structures
and the related cohomology theory from
the perspective of Lie bi-algebroids, and then address it in the context of
the holomorphic Poisson structure $\Lambda$. It quickly
sets up that the cohomology with respect to $\Lambda$
is given by the operator $\dbar_\Lambda=\adL+\dbar$ where
$\adL$ is the action of $\Lambda$ on the space of sections of the exterior product of the direct sum of
$(1,0)$-vectors and $(0,1)$-forms, and $\dbar$ is the classical operator on complex manifolds.
This frames our double complex
and its related spectral sequence.

We begin to focus on nilmanifolds in Section \ref{nilmanifolds}.
A key observation is that the Dolbeault cohomology spaces
$H^q(M, \Theta^p)$ are given by invariant elements \cite{CFP}.
It allows us to address the entire computation through that on a finite-dimensional
Lie algebra. In Section \ref{nilmanifolds}, we refine some details
regarding the action of $\adL$ on invariant elements.
This sets up the computation in Section \ref{bound}, and leads to the the proof of
Theorem \ref{phi inj}, which is also Theorem \ref{corollary theorem} stated above.

Prior to our investigation on when the map $\phi$ in Theorem \ref{corollary theorem} is
an isomorphism, in Section \ref{obstruction} we
identify the obstruction for degeneracy of the spectral sequence on the first page when the
center of the nilpotent algebra covering the manifold $M$ is of
real dimension two. Subsequently, we find that it is also sufficient for
a particular kind of holomorphic Poisson structure. Details are given in Theorem
\ref{non degenerate theorem}.

Although the degeneracy of a spectral sequence on the first page does not, in general, lead to
a  decomposition of a hypercohomology in terms of direct sums of entries in
the first page of its spectral sequence, our work leading to Theorem \ref{non degenerate theorem}
could be refined to complete a proof of Theorem \ref{Hodge},
which presents a Hodge-type decomposition.

In order to paraphrase Theorem \ref{Hodge} with precision and minimum technicality,
let $M$ be as given in the last theorem with an abelian complex
structure $J$. Let $\lieg$ be the Lie algebra of the nilpotent group covering $M$.
Let $\liec$ be the center, and $\liet$ be the quotient algebra $\lieg/\liec$. Suppose that
$\dim_\RR\liec=2$. Let $\rho$ be a non-trivial element spanning the
space of $(1,0)$-form dual to the center, i.e., $\rho\in \liec^{*(1,0)}$.
When the complex structure is abelian, then
$d\rho$ is a type $(1,1)$-form on the complexification of $\liet$ \cite{Bar, Salamon}.
Through contraction, it could be treated as a map from $\liet^{1,0}$ to $\liet^{*(0,1)}$.

\begin{theorem}\label{decomposition}
Suppose that $M$ is a nilmanifold with abelian complex structure and
suppose that $\dim_\RR\liec=2$. If the map $d\rho$ is non-degenerate,
there exists a holomorphic Poisson structure $\Lambda$ such that
the injective map $\phi$ in Theorem \ref{corollary theorem} is an isomorphism.
\end{theorem}

In order to better describe the obstruction to degeneracy of the holomorphic Poisson spectral
sequence on its first page, in Section \ref{Condition} we further analyze the necessary condition
on this obstruction. The key observations are in Proposition \ref{secondary proposition}
and Proposition \ref{main proposition}.

Given Propositions \ref{secondary proposition} and \ref{main proposition}, we examine several series
of 2-step nilmanifolds, and present examples for which the obstruction vanishes, as well as examples for
which the obstruction is non-trivial so that the map $\phi$ in Theorem \ref{corollary theorem}
fails to be an isomorphism.

\section{Complex and Generalized Complex Structures}

In this section, we review the basic background materials as seen in \cite{CGP} to set up notations.

Let $M$ be a smooth manifold. Denote its tangent bundle by $TM$ and the co-tangent bundle by $T^*M$.
A generalized complex structure on an even-dimensional manifold $M$  \cite{Marco, Hitchin-Generalized CY}
 is a
subbundle $L$ of the direct sum
${\cal T}=(TM\oplus T^*M)_{\CC}$ such that
 \begin{itemize}
 \item  $L$ and its conjugate bundle $\overline{L}$ are transversal;
 \item  $L$ is maximally isotropic with respect to the natural pairing on ${\cal T}$;
  \item the space of sections of $L$ is closed with respect to the Courant bracket.
  \end{itemize}

Given a generalized complex structure,  the pair of bundles $L$
and $\overline L$ makes a (complex) Lie bi-algebroid. The composition of the
 inclusion of $L$ and $\overline L$ in
${\cal T}$ with the natural projection onto the summand $TM_{\CC}$ becomes the anchor map
of these Lie algebroids, which we denote by $\alpha$.
Via the canonical non-degenerate pairing on the bundle $\cal T$,
 the bundle $\overline L$ is complex linearly identified to the dual of $L$.
 Therefore, the Lie algebroid differential of $\overline L$ acts on $L$, which extends to a differential on
 the exterior algebra of $L$.
 For the calculus of Lie bi-algebroids,
 we follow the conventions in \cite{Mac}. In particular,
 for any element $\Gamma$ in $C^\infty(M, \wedge^kL)$ and elements ${a}_1, \dots,
 {a}_{k+1}$ in $C^\infty(M, {\overline{L}})$, the Lie algebroid differential  of $\Gamma$  is defined by
 the Cartan formula as in exterior differential algebra, namely
 \begin{eqnarray}
&& (\dbar_{L}\Gamma)({a}_1, \dots, {a}_{k+1})
=\sum_{r=1}^{k+1}(-1)^{r+1}\alpha(a_r)(\Gamma({a}_1, \dots,{\hat{a}}_{r} ,\dots, {a}_{k+1}))\nonumber\\
&\quad& \quad +
\sum_{r<s}(-1)^{r+s}\Gamma(\lbra{a_r}{a_s}, {a}_1,
\dots,{\hat{a}}_{r} ,\dots,{\hat{a}}_{s} ,\dots,{a}_{k+1}).
\label{algebroid differential}
\end{eqnarray}
The space of sections of the bundle $\overline{L}$ is closed  with respect to the Courant bracket
if and only if $\dbar_{L}\circ \dbar_{L}=0$.

Typical examples of generalized complex structures
 are classical complex structures and symplectic structures on a manifold \cite{Marco}.
 In this section, we focus on the former.

Let $J: TM\to TM$ be an integrable complex structure on the manifold $M$.
The complexified tangent bundle $TM_\CC$ splits into the direct sum of
the bundle of $(1,0)$-vectors $TM^{1,0}$ and the bundle of $(0,1)$-vectors $TM^{0,1}$.
Their $p^{\mathrm{th}}$ exterior products are respectively denoted by $TM^{p,0}$ and $TM^{0,p}$.
Denote their dual bundles by $TM^{*(p,0)}$ and $TM^{*(0,p)}$ respectively.

Define $L= TM^{1,0}\oplus TM^{*(0,1)} $. This defines a generalized complex structure,
where its dual is its complex conjugate ${\overline{L}}=TM^{0,1}\oplus TM^{*(1,0)}$.
When one restricts the Courant bracket from the ambient bundle
${\cal T}=(TM\oplus T^*M)_{\CC}$ to the subbundles $L$ and $\overline{L}$, one then recovers the
Schouten-Nijenhuis bracket, often simply known as the Schouten bracket, in classical deformation. The
Schouten bracket between
$(1,0)$-vector fields is the Lie bracket of vector fields;
the Schouten bracket between a $(1,0)$-vector field and a $(0,1)$-form
is related to the Lie derivative of a form by a 
vector field. The Schouten bracket between two $(0,1)$-forms is equal to zero. These brackets are
extended to higher exterior products by observing the rule of exterior multiplication \cite{Mac}.

With respect to the Lie algebroid $\overline{L}$,
we get its differential $\dbar$  as defined in
(\ref{algebroid differential}),
\begin{equation}
\dbar: C^\infty(M, L ) \to
C^\infty(M, \wedge^2L).
\end{equation}
It is extended to a differential of exterior algebras:
\begin{equation}
\dbar: C^\infty(M, \wedge^pL) \to
C^\infty(M, \wedge^{p+1}L).
\end{equation}
It is an elementary exercise in computation of Lie algebroid differential
that when $\dbar$ is restricted to
$(0,1)$-forms, it is the classical $\dbar$-operator in complex manifold theory; and
\[
\dbar: C^\infty(M, TM^{*(0,1)}) \to C^\infty(M, TM^{*(0,2)})
\]
is the $(0,2)$-component of the exterior differential \cite{Brian}.
Similarly, when the Lie algebroid differential is restricted to
$(1,0)$-vector fields, then the map becomes
\[
\dbar:  C^\infty(M, TM^{1,0} )\to C^\infty(M, TM^{*(0,1)}\otimes TM^{1,0}),
\]
the {\it Cauchy-Riemann operator} as seen in \cite{Gau}.

By virtue of $L$ and $\overline{L}$ being a  Lie bi-algebroid pair, the
space $C^\infty (M, \wedge^\bullet L)$ together with the Schouten bracket,
exterior product and the Lie algebroid differential $\dbar$ form a differential Gerstenhaber
algebra  {\cite{Mac, Poon}}. In particular, if $a$ is a smooth section
of  $\wedge^{|a|}L$ 
and $b$ is a smooth section of $\wedge^{|b|}L$, then
\begin{eqnarray}
\dbar \lbra{a}{b} &=& \lbra{\dbar a}{ b}+(-1)^{|a|+1}\lbra{a}{\dbar b};\\
\dbar (a\wedge b) &=& (\dbar a)\wedge b+(-1)^{|a|}a\wedge (\dbar b).
\end{eqnarray}

Since $\dbar\circ\dbar=0$, one obtains the Dolbeault cohomology with coefficients in holomorphic
polyvector fields. Denoting the sheaf of germs of sections of
the $p^{\mathrm{th}}$ exterior power of the holomorphic tangent bundle by $\Theta^p$, we have
\[
H^\bullet(M, \wedge^\bullet \Theta)\cong
\bigoplus_{p,q\geq 0}H^q(M, \Theta^p).
\]
In subsequent computations, when $p=0$, $\Theta^p$ represents the structure
sheaf $\mathcal{O}$ of the complex manifold $M$.

Due to the compatibility between $\dbar$ and
 the Schouten bracket $\lbra{-}{-}$, and the compatibility between $\dbar$
 and the exterior product $\wedge$
 as noted above, the Schouten bracket and exterior product descend to
the cohomology space $H^\bullet(M, \Theta^\bullet)$. In other words,
the triple $(H^\bullet(M, \Theta^{\bullet,0}), \lbra{-}{-}, \wedge)$ forms a Gerstenhaber algebra.
 When we ignore the exterior product,  we call it  a Schouten algebra.
 For example, by center of the Schouten algebra, we mean the collection of  elements $A$ in
 $H^\bullet(M, \Theta^\bullet)$
  such that $\lbra{A}{B}=0$ for all $B$ in this space.

\section{Holomorphic Poisson Bi-Complex}

A holomorphic Poisson structure on a complex manifold $(M, J)$ is a
holomorphic bi-vector field $\Lambda$
such that  $\lbra{\Lambda}{\Lambda}=0$.
The corresponding bundles as a generalized complex structure for ${\Lambda}$ are the pair of
bundles of graphs $L_{\overline{\Lambda}}$ and ${\overline L}_{\Lambda}$ where
\begin{equation}
{\overline L}_{\Lambda}=\{{\overline\ell}+\Lambda(\overline{\ell}): {\overline\ell}\in {\overline L}\}.
\end{equation}

 While the pair of bundles
 $L_{\overline{\Lambda}}$ and ${\overline L}_{\Lambda}$ naturally form a Lie bi-algebroid,
 so does the pair $L$ and ${\overline L}_{\Lambda}$ \cite{LWX}. From this perspective, the Lie
 algebroid differential of the deformed generalized complex structure ${\overline L}_{\Lambda}$ acts on
 the space of sections of the bundle $L$.

Any smooth section of the bundle $L$ is the sum of a section $v$ of $TM^{1,0}$ and a section
$\oomega$ of  $TM^{*(0,1)}$.
Given a holomorphic Poisson structure $\Lambda$, define $\adL$ by
\begin{equation}
\adL(v+\oomega)=\lbra{\Lambda}{v+\oomega}.
\end{equation}

 \begin{proposition}\label{d lambda}\mbox{\rm \cite[Proposition 1]{GPR}} The action
 of the Lie algebroid differential
 of ${\overline L}_{\Lambda}$  on $L$ is given by
 \begin{equation}
 \dbar_{\Lambda}=\dbar+\adL: C^\infty(M, L)\to
 C^\infty(M, \wedge^2L).
 \end{equation}
 \end{proposition}

 The operator $\dbar_{\Lambda}$ extends to act on the exterior algebra of $TM^{1,0}\oplus TM^{*(0,1)}=L$.
  From now on, for $n\geq 0$ denote
 \begin{equation}
 K^n=C^\infty (M, \wedge^n L).
 \end{equation}
 For $n< 0$, set $K^n=\{ 0\}$.

 Since $\Lambda$ is a holomorphic Poisson structure, the closure of the space of sections of the
 corresponding bundles ${\overline L}_{\Lambda}$ is equivalent to  $\dbar_{\Lambda}\circ \dbar_{\Lambda}=0$.
 Therefore, one has a complex with $\dbar_{\Lambda}$ being a differential.
 \begin{definition}
For all $n\geq 0$, the $n^{\mathit{th}}$ holomorphic Poisson cohomology of the
holomorphic Poisson structure $\Lambda$ is the space
 \begin{equation}
 H^n_{\Lambda}(M):=\frac{{\mbox{\rm kernel of  }}\ \dbar_\Lambda: K^n \to K^{n+1} }
 {{\mbox{\rm image of }}\ \dbar_\Lambda: K^{n-1}\to K^n}.
 \end{equation}
 \end{definition}
Given Proposition \ref{d lambda}, the identity
 $\dbar_{\Lambda}\circ \dbar_{\Lambda}=0$ is equivalent to a system of three equations,
 \begin{equation}
 \dbar\circ\dbar=0, \quad
 \dbar\circ \adL+\adL\circ \dbar=0,
 \quad
 \adL\circ\adL=0.
 \end{equation}
 The first identity is equivalent to the complex structure $J$ being integrable;
 the second identity is equivalent to $\Lambda$ being holomorphic, and
 the third is equivalent to $\Lambda$ being Poisson.
Define $A^{p,q}=C^\infty(M, TM^{p,0}\otimes TM^{*(0,q)})$, then
\begin{equation}
\adL: A^{p,q}\to A^{p+1, q},
\quad
\dbar: A^{p,q} \to A^{p, q+1};
\quad
\mbox{ and }
\quad
K^n=\oplus_{p+q=n}A^{p,q}.
\end{equation}
Therefore, we obtain a bi-complex. We arrange the double indices $(p,q)$
in such a way that $p$ increases horizontally so that
$\adL$ maps from left to right, and
$q$ increases vertically so that $\dbar$ maps from bottom to top. As a result,
we obtain a {\it first quadrant } bi-complex.

\begin{definition} Given a holomorphic Poisson structure $\Lambda$, its Poisson bi-complex is
the triple $\{ A^{p,q}, \adL, \dbar\}$.
\end{definition}
It is now obvious that the (holomorphic) Poisson cohomology $H^\bullet_\Lambda(M)$
theoretically could be computed by each one of the
two naturally defined spectral sequences.
We choose a filtration given by
\[
F^pK^n=\bigoplus_{\substack{p'+q=n\\ p'\geq p}}A^{p',q}.
\]
The
lowest differential is $\dbar: A^{p,q} \to A^{p, q+1}$. Therefore, the first sheet of the
spectral sequence is the Dolbeault cohomology
\begin{equation}
E_1^{p,q}=H^q(M, \Theta^p).
\end{equation}
The first page
of the spectral sequence is given as below.
\[
\begin{array}
[c]{cccccccc}
& \cdots & \rightarrow & \cdots & \rightarrow & \cdots & \rightarrow & \\
& H^{n}(M,\mathcal{O}) & \rightarrow & H^{n}(M,\Theta) & \rightarrow
& H^{n}(M,\Theta^{2}) & \rightarrow & \\
& \cdots & \rightarrow & \cdots & \rightarrow & \cdots & \rightarrow & \\
& H^{1}(M,\mathcal{O}) & \rightarrow & H^{1}(M,\Theta) & \rightarrow
& H^{1}(M,\Theta^{2}) & \rightarrow & \\
& H^{0}(M,\mathcal{O}) & \rightarrow & H^{0}(M,\Theta) & \rightarrow
& H^{0}(M,\Theta^2) & \rightarrow &
\end{array}
\]
The differential on this page is
 \begin{equation}\label{d1}
 d_1^{p,q}=\adL: H^q(M, \Theta^p)\to H^q(M, \Theta^{p+1}).
 \end{equation}

\begin{question}\label{question on first page} When will the spectral sequence of a Poisson bi-complex
degenerate on the first page? In other words, when will
$\adL \equiv 0$ for all $p,q$?
\end{question}

It is known in the application of spectral sequences that degeneracy on the first page
does not necessarily lead to an isomorphism of the hypercohomology of degree $n$
and a direct sum of the spaces on the first page whose sum of bi-degree is equal to $n$.
A well-known case is on the degeneracy of Fr\"olicher spectral sequence on compact complex surfaces \cite{BVP}.
Therefore, the next question is also relevant.

 \begin{question}\label{question on hodge} Is there any holomorphic Poisson structure such that
 its cohomology has the following decomposition?
 \begin{equation}
 H^n_\Lambda\cong \bigoplus_{p+q=n}H^q(M, \Theta^p).
 \end{equation}
 \end{question}

 In subsequent discussion, we will use $H^{p,q}$ to represent $H^q(M, \Theta^p)$.

\section{Nilmanifolds}\label{nilmanifolds}

A compact manifold $M$ is a nilmanifold if there exists a
simply-connected nilpotent Lie group $G$
and a lattice subgroup $\Gamma$ such that $M$ is diffeomorphic to $G/\Gamma$.
We denote the Lie algebra of the group
$G$ by $\lieg$, and its center by $\liec$.
By definition,
$\lieg^0=\lieg$, and inductively, for all natural numbers $p$,
  $\lieg^{p}=\lbra{\lieg^{p-1}}{ \lieg}$. The algebra $\lieg$ is a $k$-step nilpotent algebra
  if $\lieg^{k}=\{0\}$ and $\lieg^{k-1}\neq \{0\}$. In particular,
  $\lieg^{k-1}\subseteq \liec$.
 The step of the
nilmanifold is the nilpotence of the Lie algebra $\lieg$.

A left-invariant complex
structure $J$ on $G$ is said to be {\it abelian } if on the Lie algebra $\lieg$, it satisfies the conditions
$J\circ J=-$identity and $\lbra{JA}{JB}=\lbra{A}{B}$ for all $A$ and $B$ in the Lie algebra $\lieg$ \cite{Bar,
Salamon}.
If one complexifies the algebra $\lieg$ and denotes the $+i$ and $-i$ eigen-spaces of $J$
by $\lieg^{1,0}$ and $\lieg^{0,1}$, respectively, then the invariant complex structure $J$ is abelian
if and only if the complex Lie algebra $\lieg^{1,0}$ is abelian.

Denote $\wedge^p\lieg^{1,0}$ and $\wedge^q{\lieg}^{*(0,1)}$ respectively by
$\lieg^{p,0}$ and ${\lieg}^{*(0,q)}$. We will use the notation
\[
B^{p,q}=\lieg^{p,0}\otimes {\lieg}^{*(0,q)}.
\]

On the nilmanifold $M$, we consider $\lieg^{p,0}$ as invariant $(p,0)$-vector fields and
${\lieg}^{*(0,q)}$ as invariant $(0,q)$-forms. It yields an inclusion map
\[
B^{p,q}
\hookrightarrow A^{p,q}=C^\infty(M, TM^{p,0}\otimes TM^{*(0,q)}).
\]
When the complex structure is also invariant, $\dbar$ sends
$B^{p,q}$ to
$B^{p,q+1}$.
Given an invariant complex structure and an invariant holomorphic Poisson structure $\Lambda$,
$\adL$ sends
$B^{p,q}$ to
$B^{p+1,q}$.
Restricting $\dbar$ to $B^{p,q}$, we then consider the invariant cohomology
\begin{equation}\label{hqp}
H^q( \lieg^{p,0})=
\frac
{\mbox{kernel of } \dbar: B^{p,q}\to B^{p,q+1}}
{\mbox{image of } \dbar: B^{p,q-1}\to B^{p,q}}.
\end{equation}
The inclusion map yields a homomorphism of cohomology:
\[
H^q({\lieg}^{p,0})\hookrightarrow H^q(M, \Theta^p).
\]

\begin{theorem}{\rm (See \cite{CFP, CGP})}\label{quasi isomorphic}
On a nilmanifold $M$ with an invariant
abelian complex structure, the inclusion $B^{p,q}$ in
$A^{p,q}=C^\infty(M, TM^{p,0}\otimes TM^{*(0,q)})$ induces an isomorphism of cohomology spaces. In other words,
\[
H^q({\lieg}^{p,0})\cong H^q(M, \Theta^p)=H^{p,q}.
\]
\end{theorem} 

Any element $A$ in $B^{p,q}$ acts on
$B^{p,q}$ by the Schouten bracket. We denote its action by $\ad_{A}$, i.e.,
\[
\ad_A(B)=\lbra{A}{B}.
\]
An element $A$ is said to be in the center  of the Schouten algebra
$\oplus_{p,q}B^{p,q}$ with respect to the Schouten bracket
$\lbra{-}{-}$ if and only if
$\ad_A\equiv 0$.
Similarly, an element $A$ in $H^q({\lieg}^{p,0})$ is said to be in the center of the Schouten algebra
$\oplus_{p,q}H^q({\lieg}^{p,0})$ if $\ad_A(B)$ is equal to zero on the cohomology level for
any $B$ in $\oplus_{p,q}H^q({\lieg}^{p,0})$.

For each  integer $0\leq m\leq k+1$,
define $\mathfrak{g}_{J}^{m}=\mathfrak{g}^{m}+J\mathfrak{g}^{m}.$ When the complex structure
is abelian, it is clear from various definitions that each
$\mathfrak{g}_{J}^{m}$ is a $J$-invariant ideal of $\mathfrak{g}$, and
  we have a  filtration of subalgebras:
\[
\left\{  0\right\}  =\mathfrak{g}_{J}^{k}
\subset\mathfrak{g}_J^{k-1}
\subseteq\cdots\subseteq\mathfrak{g}_{J}^{\ell+1}
\subseteq\mathfrak{g}_{J}^{\ell}
\subseteq\cdots\subseteq\mathfrak{g}_{J}^{1}\ {\subset}\ \mathfrak{g}_{J}^{0}=\mathfrak{g}.
\]
Since the center $\liec$ is $J$-invariant and it contains $\lieg^{k-1}$,
\begin{equation}\label{center}
\lie{g}^{k-1}_J\subseteq \liec.
\end{equation}

We complexify the above filtration to get
\begin{equation}\label{filtration}
\left\{  0\right\}  =\mathfrak{g}_{J,\CC}^{k}\hookrightarrow
\mathfrak{g}_{J, \CC}^{k-1}\hookrightarrow\cdots
\hookrightarrow\mathfrak{g}_{J,\CC}^{\ell+1}\hookrightarrow
\mathfrak{g}_{J,\CC}^{\ell}\hookrightarrow\cdots
\hookrightarrow\mathfrak{g}_{J,\CC}^{1}\hookrightarrow
\mathfrak{g}_{J,\CC}^{0}=\mathfrak{g}_{\CC}.
\end{equation}
With respect to the eigenspace decomposition for $J$, there exists a type decomposition for each $\ell$
\[
\mathfrak{g}_{J,\CC}^{\ell}=\mathfrak{g}_{J}^{\ell,(1,0)}\oplus\mathfrak{g}_{J}^{\ell,(0,1)}.
\]
So, the filtration (\ref{filtration}) splits into two.
One is for type $(1,0)$-vectors
\[
\left\{  0\right\}  \hookrightarrow\mathfrak{g}_{J}^{k-1,(1,0)}
\hookrightarrow\cdots\hookrightarrow\mathfrak{g}_{J}^{\ell+1,(1,0)}
\hookrightarrow\mathfrak{g}_{J}^{\ell,(1,0)}\hookrightarrow\cdots
\hookrightarrow\mathfrak{g}_{J}^{1,(1,0)}
\hookrightarrow\mathfrak{g}_{J}^{0,(1,0)}=\mathfrak{g}^{(1,0)};
\]
Another is for type $(0,1)$-vectors
\[
\left\{  0\right\}  \hookrightarrow\mathfrak{g}_{J}^{k-1,(0,1)}
\hookrightarrow\cdots\hookrightarrow\mathfrak{g}_{J}^{\ell+1,(0,1)}
\hookrightarrow\mathfrak{g}_{J}^{\ell,(0,1)}\hookrightarrow\cdots
\hookrightarrow\mathfrak{g}_{J}^{1,(0,1)}
\hookrightarrow\mathfrak{g}_{J}^{0,(0,1)}=\mathfrak{g}^{(0,1)}.
\]

\begin{lemma}\label{lemma on brackets}{\rm \cite[Lemma 4]{CFP}}
Suppose that the complex structure $J$ is abelian. Then
\begin{itemize}
\item
$\lbra{  \mathfrak{g}^{(1,0)}}{\mathfrak{g}^{(1,0)}} =0,$ and
$\lbra{ \mathfrak{g}^{(0,1)}}{\mathfrak{g}^{(0,1)}} =0.$
\item
$\lbra{  \mathfrak{g}_{J}^{a,(1,0)}}{\mathfrak{g}_{J}^{b,(0,1)}}
\subseteq\mathfrak{g}_{J, \CC}^{1+\max{\{a,b\}}}.$
\end{itemize}
In particular,
$\lbra{  \mathfrak{g}_{J}^{a,(1,0)}}{\mathfrak{g}^{(0,1)}}
\subseteq\mathfrak{g}_{J, \CC}^{a+1}$.

\end{lemma}

For the quotient space, we will adopt the notation
\begin{equation}
\mathfrak{t}^{(1,0)}_{\ell}
=\lie{g}_{J}^{\ell-1,(1,0)}/
\lie{g}_{J}^{\ell,(1,0)}.
\end{equation}
Choose a vector space isomorphism so that the short exact sequence of Lie algebras
\[
0\to \lie{g}_{J}^{\ell, (1,0)}\to \lie{g}_{J}^{\ell-1, (1,0)} \to \lie{g}_{J}^{\ell-1,(1,0)}/
\lie{g}_{J}^{\ell,(1,0)}\to 0
\]
is turned into a direct sum of vector spaces.
\[
\mathfrak{g}_{J}^{\ell-1,(1,0)}\cong\mathfrak{t}_{\ell}^{(1,0)}
\oplus\mathfrak{g}_{J}^{\ell,(1,0)}.
\]
Inductively,
\begin{equation}\label{g10}
\mathfrak{g}^{(1,0)}=\mathfrak{t}_{1}^{(1,0)}
\oplus\mathfrak{t}_{2}^{(1,0)}\oplus\cdots\oplus\mathfrak{t}_{k}^{(1,0)}.
\end{equation}
Similarly,
\[
\mathfrak{g}^{(0,1)}=\mathfrak{t}_{1}^{(0,1)}
\oplus\mathfrak{t}_{2}^{(0,1)}\oplus\cdots\oplus\mathfrak{t}_{k}^{(0,1)}.
\]
 We remark that $\mathfrak{t}_{k}^{(1,0)}=\mathfrak{g}_J^{k-1,(1,0)}$. Let
 $\lie{t}_\ell^{*(0,1)}$ be the dual space of $\lie{t}_{\ell}^{(0,1)}$, and
 $\lie{t}_{\ell}^{*(1,0)}$ the dual space of $\lie{t}_{\ell}^{(1,0)}$.

\begin{lemma}\label{holomorphic}{\rm \cite[Proposition 1]{CFP}}
When the complex structure is abelian and
$\lieg$ is $k$-step nilpotent, then
 \begin{enumerate}
 \item $\overline{\partial}\lieg^{*(0,1)}=\{ 0\}$.
 \item $\overline{\partial}\mathfrak{t}_{k}^{1,0}=\{0\}$.
 \item $\overline{\partial}\mathfrak{t}_{k-1}^{1,0}\subseteq
\oplus_{a\leq k-1}(\lie{t}_{k}^{1,0}\otimes \lie{t}_{a}^{*(0,1)}).$
 \item For all
  $1\leq \ell \leq k-2$,
$
\overline{\partial}\mathfrak{t}_{\ell}^{1,0}\subseteq
\bigoplus_{a\leq \ell}(\lie{t}_{\ell+1}^{1,0}\otimes \lie{t}_{a}^{*(0,1)})
\,{\large\oplus}\, \bigoplus_{a> \ell}\left( \lie{t}_{a+1}^{1,0}\otimes \lie{t}_{a}^{*(0,1)}
\right).$
\end{enumerate}
\end{lemma}

As a consequence, all $(0,q)$-forms are $\dbar$-closed.
\begin{lemma}\label{0q forms} There is an isomorphism
\[
H^q(\lieg^{0,0})=B^{0,q}=\wedge^q\lieg^{*(0,1)}.
\]
\end{lemma}

\section{A Bound on Holomorphic Poisson Cohomology}\label{bound}

Due to Theorem \ref{quasi isomorphic}, our investigation on the holomorphic Poisson spectral sequence
for any invariant holomorphic Poisson structure $\Lambda$
 is reduced to an analysis on the adjoint action of $\Lambda$ on the invariant cohomology
 $H^q(\lieg^{p,0})$.

 Given the bi-degree, we next create a filtration of cohomology.
For any natural number $n$, $B^n=\oplus_{p+q=n}B^{p,q}$. Define
\begin{equation}
F^mB^n=\bigoplus_{\substack{m\leq j\leq n\\ j+\ell=n}}B^{j,\ell}
=B^{n, 0}\oplus B^{n-1, 1}\oplus \cdots \oplus B^{m+1, n-m-1}\oplus B^{m, n-m}.
\end{equation}
It follows that $F^0B^n=B^n$, $F^{n}B^n=B^{n,0}$. Define $F^mB^n=\{ 0\}$ when $m\geq n+1$.
Therefore, we have a filtration
\begin{equation}
B^n=F^0B^n \supset F^1B^n \supset \cdots \supset F^kB^n
\supset \cdots \supset F^nB^n \supset F^{n+1}B^n=\{0\}.
\end{equation}
Define
\begin{equation}
F^mZ^n=\ker \dbar_\Lambda: F^mB^n\to B^{n+1},
\quad
F^mC^n=\image \dbar_\Lambda: F^mB^{n-1}\to B^{n}.
\end{equation}
As $\dbar_\Lambda$ maps $B^{p,q}$ to $B^{p+1, q}\oplus B^{p,q+1}$, it is apparent that
$\dbar_\Lambda$ maps $F^mB^{n-1}$ to $F^mB^{n}$. Therefore, $F^mZ^n\supseteq F^mC^n$, and
the quotient below is well defined
\begin{equation}
F^mH^n=\frac{F^mZ^n}{F^mC^n}.
\end{equation}
As $F^{m-1}Z^n\cap F^mC^n=F^{m-1}C^n$, we have an inclusion
$F^{m-1}H^n\supset F^mH^n$, and hence a filtration
\begin{equation}
H^n_{\Lambda}=F^0H^n \supset F^1H^n \supset \cdots
\supset F^kH^n \supset \cdots \supset F^nH^n \supset F^{n+1}H^n=\{0\}.
\end{equation}
This induces a vector space isomorphism
\begin{equation}
H^n_{\Lambda}\cong \bigoplus_{m=0}^n\frac{F^mH^n}{F^{m+1}H^n}.
\end{equation}
Each element in $\frac{F^mH^n}{F^{m+1}H^n}$ is represented by an element $\alpha$ in $F^mB^n$ such that
$\dbar_\Lambda\alpha=0$. Let us expand $\alpha$ as the sum of elements in $B^{p,q}$. It is given by
\begin{equation}
\alpha=\sum_{\substack{m\leq j\leq n\\ j+\ell=n}}\alpha^{j,\ell}
=\alpha^{n, 0}+\alpha^{n-1, 1}+\cdots +\alpha^{m+1, n-m-1}+\alpha^{m,n-m}.
\end{equation}
Therefore,
\begin{eqnarray*}
\dbar_\Lambda\alpha&=&\adL\alpha+\dbar\alpha \\
&=&\adL \alpha^{n, 0}+\adL\alpha^{n-1, 1}+\cdots +\adL\alpha^{m+1, n-m-1}+\adL\alpha^{m,n-m}\\
&\quad & \quad \hspace{0.5in}
\dbar\alpha^{n, 0}+\dbar\alpha^{n-1, 1}+\cdots +\dbar\alpha^{m+1, n-m-1}+\dbar\alpha^{m,n-m}.
\end{eqnarray*}
Since the complex structure is abelian, $\adL\alpha^{n,0}=0$. By collecting terms according to the
bi-degrees in $B^{p,q}$, we obtain
\[
\dbar_{\Lambda}\alpha=(\adL\alpha^{n-1, 1}+\dbar\alpha^{n, 0}) +
\cdots +(\adL\alpha^{m,n-m}+\dbar\alpha^{m+1, n-m-1})+\dbar\alpha^{m,n-m}.
\]
Therefore, $\dbar_{\Lambda}\alpha=0$ if and only if
\[
\adL\alpha^{n-1, 1}+\dbar\alpha^{n, 0}=0, \dots, \adL\alpha^{m,n-m}+\dbar\alpha^{m+1, n-m-1}=0
,\quad {\mbox{and}}\quad  \dbar\alpha^{m,n-m}=0.
\]
In particular, $\alpha^{m, n-m}$ represents an element in the cohomology space
$H^{m, n-m}=H^{n-m}(\lieg^{m,0}).$

We use the notation $[-]_\Lambda$ to indicate equivalence classes with respect to the total
cohomology operator $\dbar_{\Lambda}=\adL+\dbar$, and use $[-]_{\dbar}$ to
indicate equivalence classes with respect to $\dbar$. The observation in the last paragraph indicates that
for any $n\geq0$ and $0\leq m\leq n$, we could define a map
\begin{equation}\label{phi map}
\phi_{m, n-m}:F^mH^n\to H^{m, n-m}, \quad {\rm{where} } \quad \phi_{m, n-m}([\alpha]_\Lambda)
=[\alpha^{m,n-m}]_{\dbar}.
\end{equation}

To check that this map is well defined, consider $\alpha$ and $\beta$ in $F^mB^n$ such that
$[\alpha]_\Lambda=\classl{\beta}$. As their difference is in $F^mC^n$, there exists
$\gamma\in F^m B^{n-1}$ such that $\alpha=\beta+\dbar_{\Lambda}\gamma$. As
\[
\gamma=
\sum_{\substack{ m\leq j\leq n-1 \\ j+\ell=n-1}}
\gamma^{j,\ell}
=\gamma^{n-1, 0}+\cdots +\gamma^{m+1, n-1-m-1}+\gamma^{m,n-1-m},
\]
and $\dbar_{\Lambda}\gamma^{j,\ell}=\adL\gamma^{j,\ell}+\dbar\gamma^{j,\ell}$,
the only term of $\dbar_{\Lambda}\gamma$ in $B^{m, n-m}$
is given by $\dbar\gamma^{m, n-1-m}$.
It follows that
\[
\alpha^{m,n-m}=\beta^{m,n-m}+\dbar\gamma^{m, n-1-m},
\]
and hence $\classd{\alpha^{m,n-m}}=\classd{\beta^{m,n-m}}$. In other words, the map
$\phi_{m, n-m}$ is well defined at cohomology level.

Furthermore, suppose that $\classl{\alpha}$ is in the kernel of the map $\phi_{m, n-m}$. Then there
exists $\gamma^{m, n-m-1}\in B^{m, n-m-1}$ such that $\alpha^{m, n-m}=\dbar\gamma^{m, n-m-1}$, so
\[
\alpha=
\alpha^{n, 0}+\alpha^{n-1, 1}+\cdots +\alpha^{m+1, n-m-1}+\dbar\gamma^{m, n-m-1}
\]
and $\dbar_{\Lambda}\alpha=0$. As $\classl{\alpha}=\classl{\alpha-\dbar_\Lambda\gamma^{m, n-m-1}}$ and
\[
\alpha-\dbar_\Lambda\gamma^{m, n-m-1}
=\alpha^{n, 0}+\alpha^{n-1, 1}+\cdots +\alpha^{m+1, n-m-1}-\adL\gamma^{m, n-m-1}
\]
is contained in $F^{m+1}H^n$, the kernel of $\phi_{m, n-m}$ is contained in
$F^{m+1}H^n$. On the other hand, by inspecting the bi-degree, $F^{m+1}H^n$ is contained in
the kernel of $\phi_{m, n-m}$, so  $F^{m+1}H^n$ is equal to the kernel of the map $\phi_{m, n-m}$.

For each $n$, by taking the direct sum of linear maps $\phi=(\phi_{n,0}, \dots, \phi_{m, n-m})$,
we obtain the next observation.

\begin{theorem}\label{phi inj}
For any invariant holomorphic Poisson structure $\Lambda$ on a nilmanifold with
an abelian complex structure, there exists an injective map $\phi$ from $H^n_\Lambda$ into
 $\oplus_{p+q=n}H^{p,q}=\oplus_{p+q=n}H^q(\lieg^{p,0}).$
 In particular,
 \[
 \dim H^n_\Lambda\leq \sum_{p+q=n}\dim H^{p,q}=\sum_{p+q=n}\dim H^q(\lieg^{p,0}).
 \]
\end{theorem}

Now an obvious question is whether the map $\phi$
 is surjective. On a nilmanifold with abelian complex
structure, there is an obvious set of cohomology spaces as seen in Lemma \ref{0q forms}, namely
$H^{0,q}=H^q(\lieg^{0,0})=B^{0,q}$.
Suppose that $\alpha^{0,1}$ is in $B^{0,1}$. If the map $\phi$ is surjective,
then there exists $\alpha^{1,0}\in B^{1,0}$ such that
$\adL(\alpha^{1,0}+\alpha^{0,1})=0$. It means
\[
\adL(\alpha^{1,0})+\left(\dbar\alpha^{1,0} +\adL(\alpha^{0,1})\right)+\dbar\alpha^{0,1}=0.
\]
Since the complex structure is abelian, $\adL\alpha^{1,0}=0$ and $\dbar\alpha^{0,1}=0$. It follows that
\[
\dbar\alpha^{1,0} +\adL\alpha^{0,1}=0.
\]
In other words, $\adL\alpha^{0,1}$ is $\dbar$-exact, a non-trivial criterion for the holomorphic
Poisson spectral sequence degenerates on its first page.

\section{A Class of Holomorphic Poisson Structures}
Since $\dbar_\Lambda=\adL+\dbar$, if the action of $\adL$ on $B^{p,q}$ is equal to zero for all $p,q$,
the action of $\dbar_\Lambda$ on $B^{p,q}$ is equal to the action of $\dbar$. In such case,
$\phi$ in Theorem \ref{phi inj} is apparently an isomorphism.
This would be the case if $\Lambda$ is a non-trivial element in
 $\liec^{2,0}=\wedge^2 \liec^{1,0}$.

From now on, we consider the case when $\dim_{\CC}\liec^{1,0}=1$.
In terms of the notation in the
decomposition (\ref{g10}), let $\liet_k^{1,0}=\liec^{1,0}$.
We will use the notation $\liet^{1,0}=\liet_1^{1,0}\oplus \cdots \oplus \liet_{k-1}^{1,0}$.

\begin{lemma}\label{dbar closed} Suppose that $V\in \liet_{k}^{1,0}$ and $T\in\liet_{k-1}^{1,0}$. If
$\dim_{\CC}\liec=1$, then $\Lambda= V\wedge T$ is a holomorphic Poisson structure.
\end{lemma}
\bproof By item 2 of Lemma \ref{holomorphic}, $\dbar(V\wedge T)=(\dbar V)\wedge T-V\wedge (\dbar T)
=-V\wedge (\dbar T)$. By item 3 of the same lemma $V\wedge (\dbar T)$ is an element in
$\oplus_{a\leq k-1}(\lie{t}_{k}^{2,0}\otimes \lie{t}_{a}^{*(0,1)})$.
Since $\liet_k^{1,0}=\liec^{1,0}$ when $\dim\liec^{1,0}=1$,
 $\liet_k^{2,0}=\{0\}.$ Therefore, $V\wedge (\dbar T)$ vanishes.

Since the complex structure is abelian $\lbra{V\wedge T}{V\wedge T}=0$,
$\Lambda=V\wedge T$ is an invariant holomorphic Poisson structure.
\eproof

Next, we refine an observation in \cite{Console-Fino, Salamon}.
\begin{lemma}\label{Schouten}
The Schouten bracket between elements in $\liet_{a}^{1,0}$ and $\liet_h^{*(0,1)}$ are
given below.
\begin{itemize}
\item $\lbra{\liet_a^{1,0}}{\liet_h^{*(0,1)}}=\{ 0\}$,  if $h\leq a$.
\item $\lbra{\liet_a^{1,0}}{\liet_h^{*(0,1)}}\subset \oplus_{b<h}\liet_b^{*(0,1)}$,  if $h= a+1$.
\item  $\lbra{\liet_a^{1,0}}{\liet_h^{*(0,1)}}\subset \liet_{h-1}^{*(0,1)}$,  if $h\geq a+2$.
\end{itemize}
\end{lemma}
\bproof
Suppose that $V_a\in \liet_a^{1,0}$, ${\overline V}_b\in \liet_b^{0,1}$ and
$\oomega^h\in \liet_h^{*(0,1)}$, then
\[
\lbra{V_a}{\oomega^h}({\overline V}_b)=(\iota_{V_a}d\oomega^h){\overline V}_b=-
\oomega^h(\lbra{V_a}{{\overline V}_b}).
\]
By Lemma \ref{lemma on brackets}, $\lbra{V_a}{{\overline V}_b}$ is contained in
$\liet_{\max\{a, b\}+1}^{1,0}$. Therefore, the evaluation
$\oomega^h(\lbra{V_a}{{\overline V}_b})$ is non-zero only if
$h=\max\{a, b\}+1$. The first bullet point is now evident. Furthermore,
when $h=a+1$, then the evaluation may be non-zero only if $b \leq a$, i.e, $b< h$.
It yields the second observation.
When $h\geq a+2$, the evaluation may be non-zero only if $\max\{a, b\}+1=b+1$ and hence
$b=h-1$.    It yields the last observation of this lemma.
\eproof

\begin{lemma}\label{top wedge}{\rm \cite[Corollary 1]{CFP}} Suppose that
$\lieg$ is $k$-step nilpotent, then
\begin{enumerate}
\item For all  $m$,
$\lbra{\liec^{1,0}}{\lieg^{*(0,1)}} =\{0\}.$
\item $\lbra{\mathfrak{t}_{k-1}^{1,0}}{\liet^{*(0,1)}}=\{0\}.$
\item $\lbra{\mathfrak{t}_{k-1}^{1,0}}{\liec^{*(0,1)}} \subseteq
\liet^{*(0,1)}.$
\item Let $V\in \liet_k^{1,0}$ and $T\in \liet_{k-1}^{1,0}$, and $\Lambda=V\wedge T$, then
$\adL(\liet^{*(0,1)})=\{0\}$; and
$\adL(\liec^{*(0,1)})\subseteq \liec^{1,0}\otimes\liet^{*(0,1)}$.
\end{enumerate}
\end{lemma}
\bproof With $\liec^{1,0}=\liet_k^{1,0}$ and $\liet^{*(0,1)}=\bigoplus_{1\leq \ell\leq k-1}\liet_{\ell}^{*(0,1)}$,
the first three items are the direct consequence of Lemma \ref{Schouten}.
To address the last item, we note that the first item implies that for all $\ell$,
\[
\adL(\liet_\ell^{*(0,1)})=V\wedge \ad_T(\liet_\ell^{*(0,1)}).
\]
By the first item, the only non-trivial case is when $\ell=k$. The result then follows from
the second item.
\eproof

\section{Obstruction for Degeneracy}\label{obstruction}

Suppose that $\Lambda=V\wedge T$ as given in the last section, so
it is an invariant holomorphic Poisson structure.
We will denote the dual element of $V$ by $\rho$, which then spans $\liet_k^{*(1,0)}=\liec^{*(1,0)}$.

Prior to our investigation on the map $\phi$ in Theorem \ref{phi inj},
we explore the necessary condition for
the spectral sequence of the bi-complex of
$\Lambda$ to degenerate on the first page.

Consider the action of $\adL$ on the $0^{\mathrm{th}}$ column on the first page
of the spectral sequence.
\[
\adL: H^q(\lieg^{0,0})\rightarrow  H^q(\lieg^{1,0}).
\]
By Lemma \ref{0q forms}, $H^q(\lieg^{0,0})=B^{0,q}=
\wedge^{q}(\oplus_{1\leq \ell\leq k}\liet_\ell^{*(0,1)})$.
By Lemma \ref{Schouten}, the action of $\adL$ on $B^{0,q}$
is equal to zero, except possibly when an element
is of the form ${\overline{\rho}}\wedge{\overline{\Omega}}$ where $\Oomega$ is in
$\wedge^{q-1}(\oplus_{1\leq \ell\leq k-1}\liet_\ell^{*(0,1)})$.

In particular,  $\overline{\rho}$  represents an element
in $H^1(\lieg^{0,0})$. By Lemma \ref{top wedge},
\[
\adL\overline{\rho}=\lbra{V\wedge T}{\overline\rho}=V\wedge \lbra{T}{\overline\rho}
\in \liec^{1,0}\otimes\liet^{*(0,1)}
\]
It is
$\dbar$-exact only if there exists
$X$ in $\liet^{1,0}$ such that
\begin{equation}\label{dbar X definition}
V\wedge \lbra{T}{\overline\rho}=
\lbra{\Lambda}{\overline{\rho}} =\dbar X.
\end{equation}

To summarize our computation so far, we note the following.
\begin{lemma} If the holomorphic Poisson spectral sequence for $\Lambda=V\wedge T$ degenerates on
the first page, there exists a vector $X$ in $\liet^{1,0}$ such that
$\adL{\overline\rho}=\dbar X$.
\end{lemma}

We will examine the conditions under which such a vector $X$ exists. We demonstrate next that its
existence is the only obstruction for the holomorphic Poisson spectral sequence to degenerate on
its first page.

Given the dimension constraint on $\liec^{1,0}$, Lemma \ref{top wedge} shows that the action of
$\ad_{\Lambda}$ is possibly non-trivial only when it acts on components of type
\[
\Upsilon \in \lieg^{p,0}
\otimes\liet^{*(0,q-1)}\otimes \liec^{*(0,1)}.
\]
Given such a $\Upsilon$, there exist finitely many
$\overline{\Omega}_i$ in $\liet^{*(0,q-1)}$
and the same number of $\Theta_i$ in
$\lieg^{p,0}$ such that
\[
\Upsilon =\overline{\rho}\wedge \sum_{i}(\Theta_i\wedge\overline{\Omega}_i).
\]
Therefore, $\adL{\Upsilon}$ is equal to
\[
\lbra{\Lambda}{\overline{\rho}\wedge \sum_i(\Theta_i\wedge\overline{\Omega}_i)}
=V \wedge \lbra{T}{\overline{\rho}\wedge \sum_i(\Theta_i\wedge\overline{\Omega}_i)}
-T \wedge \lbra{V}{\overline{\rho}\wedge \sum_i(\Theta_i\wedge\overline{\Omega}_i)}.
\]

By Lemma \ref{holomorphic} and Lemma \ref{top wedge},
the last term on the right hand side is identically zero. By Lemma
\ref{top wedge} again,  the first term on the right hand side is equal to
\[
V\wedge \lbra{T}{\overline{\rho}}\wedge \sum_i(\Theta_i\wedge\overline{\Omega}_i)
= \lbra{V\wedge T}{\overline{\rho}}\wedge \sum_i(\Theta_i\wedge\overline{\Omega}_i)
\]
If  there exists $X$ in $\liet^{1,0}$
such that $\lbra{V\wedge T}{\overline{\rho}}=\lbra{\Lambda}{\overline{\rho}}=\dbar X$,
\begin{equation}\label{dbar closed 1}
\lbra{\Lambda}{\Upsilon}=
\lbra{V\wedge T}{\overline{\rho}}\wedge \sum_i(\Theta_i\wedge\overline{\Omega}_i)
=\dbar X \wedge \left(\sum_i(\Theta_i\wedge\overline{\Omega}_i)\right).
\end{equation}
Given that $\overline\rho$ is $\dbar$-closed and every element in $B^{0,q-1}$ is $\dbar$-closed,
\begin{equation}\label{dbar closed 2}
\dbar\Upsilon=\dbar(\overline{\rho}\wedge  \sum_i(\Theta_i\wedge\overline{\Omega}_i)   )
=-\overline{\rho}\wedge \dbar( \sum_i(\Theta_i\wedge\overline{\Omega}_i) )
=-\overline{\rho}\wedge ( \sum_i(\dbar\Theta_i)\wedge\overline{\Omega}_i )
\end{equation}
By Lemma \ref{holomorphic}, $\dbar\Theta_i$ does not have any $\overline{\rho}$ component.
Therefore, $\dbar\Upsilon=0$ if and only if
$\dbar( \sum_i(\Theta_i\wedge\overline{\Omega}_i) )=0$.
It follows that Identity (\ref{dbar closed 1}) is further transformed to
\begin{equation}
\lbra{\Lambda}{\Upsilon}=
\dbar (X \wedge \sum_i(\Theta_i\wedge\overline{\Omega}_i) ).
\end{equation}

\begin{lemma}\label{dbar exact}
 When $\adL{\overline\rho}=\dbar X$, $\adL$ is a zero map on $H^{p,q}$ for
all $p,q$. Moreover, if
$\Upsilon =\overline{\rho}\wedge \sum_{j}(\Theta_j\wedge\overline{\Omega}_j)$ is $\dbar$-closed,
then $\adL(\Upsilon)$ is $\dbar$-exact. To be precise,
\[
\adL(\Upsilon)=\dbar (X \wedge \sum_i(\Theta_i\wedge\overline{\Omega}_i) ).
\]
\end{lemma}

As a consequence, the necessary condition for the spectral sequence of
the bi-complex of $\Lambda$ to degenerate on the first page is also sufficient.

\begin{theorem}\label{non degenerate theorem} Let $M=G/\Gamma$ be a
k-step nilmanifold with an abelian complex structure.
Let $\liec$ be the center of the Lie algebra $\lieg$ of the simply connected Lie group $G$.
Let $\lieg^{1,0}=\liet^{1,0}\oplus \liec^{1,0}$ be the space of invariant $(1,0)$-vectors.
Assume that $\dim\liec^{1,0}=1$.
Suppose that $\Lambda=V\wedge T$ where $V$ is in $\liec^{1,0}$ and $T$ is in
$\liet_{k-1}^{1,0}$.
Let $\overline\rho$ span $\liec^{*(0,1)}$. The spectral
sequence of the bi-complex of $\Lambda$ degenerates on the first page if and only if
 $\ad_{\Lambda}\overline\rho$ is $\dbar$-exact.
\end{theorem}

This theorem answers Question \ref{question on first page} for nilmanifolds with an abelian complex
structure.

\section{Hodge-Type Decomposition}

Given the discussion on degeneracy on first page, the answer to Question \ref{question on hodge} becomes
useful. We now examine when the map $\phi_{m, n-m}$ is surjective. Suppose that
$\alpha^{m, n-m}$ represents a class in $H^{m, n-m}$. In particular, it is in
\[
B^{m, n-m}=\lieg^{m,0}\otimes \wedge^{n-m}( \liet^{*(0,1)}\oplus\liec^{*(0,1)})
= \lieg^{m,0}\otimes (\liet^{*(0, n-m)}\,\oplus\,\liec^{*(0,1)}\otimes \liet^{*(0,n-m-1)}).
\]

There exists finitely many $\Theta_i$ and $\Pi_j$ in $\lieg^{m,0}$, $\Oomega_i$ in
$\liet^{*(0,n-m-1)}$ and ${\overline\Gamma}_j$ in $\liet^{*(0,n-m)}$ such that
\begin{equation}
\alpha^{m, n-m}=\sum_j\Pi_j\wedge {\overline\Gamma}_j
+{\overline\rho}\wedge \left(\sum_{i}\Theta_i\wedge \Oomega_i\right) .
\end{equation}
By Lemma \ref{holomorphic},
\begin{eqnarray*}
\dbar\alpha^{m, n-m}
&=&\dbar\left(\sum_j\Pi_j\wedge {\overline\Gamma}_j\right)
-{\overline\rho}\wedge \dbar\left(\sum_i\Theta_i\wedge \Oomega_i \right)
\\
&=&\sum_j(\dbar\Pi_j)\wedge {\overline\Gamma}_j
-{\overline\rho}\wedge \left(\sum_i(\dbar\Theta_i)\wedge \Oomega_i\right).
\end{eqnarray*}
By the same lemma,
\begin{eqnarray*}
{\overline\rho}\wedge \left(\sum_i(\dbar\Theta_i)\wedge \Oomega_i\right)  &\in&
\lieg^{m,0}\otimes \liec^{*(0,1)}\otimes \liet^{*(0,n-m)}\\
(\dbar\Pi_i)\wedge {\overline\Gamma}_j&\in &
\lieg^{m,0}\otimes \liet^{*(0, n-m+1)}.
\end{eqnarray*}
As they are in different components, when $\dbar\alpha^{m, n-m}=0$, each of them is equal to zero.
\begin{equation}\label{dbar closed 3}
\dbar(\sum_j\Pi_j\wedge {\overline\Gamma}_j)=0, \quad
 {\overline\rho}\wedge\left(\sum_i(\dbar\Theta_i)\wedge \Oomega_i\right)=0.
\end{equation}
Therefore, to prove that the map $\phi_{m, n-m}$ is surjective, it suffices to examine two cases
independently; that is, either when
\[
\alpha^{m,n-m}=\sum_j\Pi_j\wedge {\overline\Gamma}_j\in
\lieg^{m,0}\otimes  \liet^{*(0,n-m)},
\]
or when
\[
\alpha^{m,n-m}={\overline\rho}\wedge \left(\sum_i\Theta_i\wedge \Oomega_i\right)
\in
\lieg^{m,0}\otimes \liec^{*(0,1)}\otimes \liet^{*(0,n-m-1)}.
\]

In the former case, Lemma \ref{top wedge} implies
$\adL\alpha^{m, n-m}=0$, and hence $\dbar_\Lambda\alpha^{m, n-m}=0$. Therefore,
$\alpha^{m, n-m}$ represents a class in $F^mH^n$ such that
$\phi_{m, n-m}(\classl{\alpha^{m, n-m}})=\classd{\alpha^{m, n-m}}$.

In the latter case, by Identity (\ref{dbar closed 3}), $\dbar\alpha^{m, n-m}=0$ implies that
$\sum_i\Theta_i\wedge\Oomega_i$ is $\dbar$-closed. Hence,
\begin{eqnarray*}
\adL\alpha^{m, n-m}&=& \adL\left({\overline\rho}\wedge \sum_i(\Theta_i\wedge\Oomega_i)\right)
=(\adL{\overline\rho})\wedge \left(\sum_i(\Theta_i\wedge\Oomega_i)\right)\\
&=&(\dbar X)\wedge\left(\sum_i(\Theta_i\wedge\Oomega_i)\right)
= \dbar \left( X\wedge \sum_i(\Theta_i\wedge\Oomega_i)\right).
\end{eqnarray*}
Define
\[
\alpha^{m+1, n-(m+1)}=-X\wedge \sum_i(\Theta_i\wedge\Oomega_i).
\]
Since it is in $\lieg^{m+1,0}\otimes \liet^{*(0,n-m-1)}$, Lemma \ref{top wedge} implies that
$\adL(\alpha^{m+1, n-(m+1)})=0$. Therefore,
\begin{eqnarray*}
&& \dbar_{\Lambda}(\alpha^{m+1, n-(m+1)}+\alpha^{m, n-m})\\
&=& \adL(\alpha^{m+1, n-(m+1)})+\dbar(\alpha^{m+1, n-(m+1)})
+\adL(\alpha^{m, n-m})+\dbar(\adL(\alpha^{m, n-m}))\\
&=&-\dbar\left( X\wedge \sum_i(\Theta_i\wedge\Oomega_i)\right)
+\dbar \left( X\wedge \sum_i(\Theta_i\wedge\Oomega_i)\right)=0.
\end{eqnarray*}
Finally,
\[
\alpha^{m+1, n-(m+1)}+\alpha^{m, n-m}=-X\wedge \sum_i(\Theta_i\wedge\Oomega_i)
+
{\overline\rho}\wedge \left(\sum_i\Theta_i\wedge \Oomega_i\right)
\]
represents a class in $F^mH^n$ such that
$\phi_{m, n-m}(\classl{\alpha^{m+1, n-(m+1)}+\alpha^{m, n-m}})
=\classd{\alpha^{m, n-m}}$,
and the map $\phi_{m, n-m}$ is surjective.

\begin{theorem}\label{Hodge}
Let $M=G/\Gamma$ be a $k$-step nilmanifold with an invariant abelian complex
structure. Let $\lieg$ be the Lie algebra of the covering group,
 $\liec$ the center.  Let $V$ be in $\liec^{1,0}$ and $\rho$ be the dual of $V$.
 Let $T$ be in $\liet_{k-1}^{1,0}$. If $\adL({\overline{\rho}})$ is
$\dbar$-exact, the holomorphic Poisson cohomology of $\Lambda=V\wedge T$ has a Hodge-type
decomposition:
\[
H^n_{\Lambda}\cong \oplus_{p+q=n}H^q(M, \Theta^p).
\]
\end{theorem}

\section{$\dbar$-Exactness of $\adL\overline\rho$}\label{Condition}
In this section, we explore when there exists $X$ such that
\begin{equation}\label{exactness 1}
\lbra{V\wedge T}{\overline{\rho}}=\adL\overline{\rho}=\dbar  X.
\end{equation}
On the left of the equality above,
\begin{equation}\label{exactness 2}
\lbra{V\wedge T}{\overline{\rho}}
=V\wedge  \lbra{ T}{\overline{\rho}}=V\wedge \iota_{T}d{\overline{\rho}}.
\end{equation}
To compute $\dbar X$, we applies the Cartan formula to evaluate $\dbar X$ on a generic element
$\overline{Y}$ in $\liet^{0,1}$ and on $\rho$.
\begin{eqnarray*}
\dbar X(\rho, {\overline{Y}}) &=& -X(\lbra{\rho}{\overline{Y}})
=X(\lbra{\overline{Y}}{\rho})=X(\iota_{\overline{Y}}d\rho)\\
&=&d\rho( \overline{Y}, X)=-d\rho(X, \overline{Y})=-(\iota_Xd\rho)(\overline{Y}).
\end{eqnarray*}
Therefore, $\dbar X=-V\wedge\iota_Xd\rho$. Now comparing (\ref{exactness 1}) with
(\ref{exactness 2}), we obtain a rather simple identity below.
\begin{equation}\label{dbar exact equation}
\iota_{T}d{\overline{\rho}}=-\iota_X d\rho.
\end{equation}

Since the complex structure $J$ is abelian, $d\rho$ is a type $(1,1)$-form.
So is $d\overline{\rho}$. We could treat their contractions with elements in
$\liet^{1,0}$  as linear maps.
\begin{equation}
d\rho, \quad d{\overline\rho} \quad : \liet^{1,0}\to \liet^{*(0,1)}.
\end{equation}

Suppose that $d\overline\rho$ has a non-trivial kernel and
$\iota_{T}d{\overline{\rho}}=0$. In such case,
$\ad_{V\wedge T}\overline\rho=0$. However, since $d\overline\omega^j=0$ for all $j$,
$\ad_{V\wedge T}\overline\omega^j=0$. As the complex structure is abelian, the adjoint action of
$V\wedge T$ on $\lieg^{1,0}$ is identically zero. Therefore, $\ad_{V\wedge T}$ is identically zero
on $B^{p,q}$ for all $p, q\geq 0$.  In such a case, the action of
 $\ad_{\Lambda}$ is trivial. Hence as an action on $B^{p,q}$, $\dbar_{\Lambda}=
 \dbar$.
  The Poisson cohomology is simply the direct sum of Dolbeault cohomology,
\[
H^n_{\Lambda}(M)\cong \bigoplus_{p+q=n} H^q(\lieg^{p,0})=\bigoplus_{p+q=n} H^q(M, \Theta^p).
\]

On the other hand,  $d\rho(\liet^{1,0})$  is a proper subspace of
$\liet^{*(0,1)}$ if $d\rho$ degenerates. If $T$ is not in the kernel of  $d{\overline\rho}$ and
if $\iota_Td{\overline\rho}$ is in the complement of $d\rho(\liet^{1,0})$,
Equation (\ref{dbar exact equation}) does not  have a solution. In the next section, we will
present an example to demonstrate that such situation does occur.
We summarize our discussion when $d\rho$ degenerates as below.

\begin{proposition}\label{secondary proposition}
If $d\overline\rho$ degenerates and $T$ is in its kernel, then $\dbar_{\Lambda}=\dbar$
for $\Lambda=V\wedge T$. In such a case, the holomorphic Poisson cohomology has a
Hodge-type composition.
\end{proposition}

If $d{\overline\rho}$ is non-degenerate, then so is $d\rho$. Therefore, for any $T$ the equation
 (\ref{dbar exact equation}) always has a unique solution.

\begin{proposition}\label{main proposition}
Let $M=G/\Gamma$ be a $k$-step nilmanifold with an abelian complex structure.
Let $\liec$ be the center of the Lie algebra $\lieg$ of the simply connected Lie group $G$.
Assume that $\dim\liec^{1,0}=1$. Let $V$ span $\liec^{1,0}$
and $\rho$ span the dual space.
If $d\rho$ is non-degenerate, then for any $T$ in $\liet_{k-1}^{1,0}$ and
 $\Lambda=V\wedge T$, $\adL(\overline{\rho})$ is $\dbar$-exact.
 In particular, the holomorphic Poisson cohomology for such $\Lambda$
 has a Hodge-type decomposition.
 \end{proposition}

\section{Examples}\label{examples}

We now focus on 2-step nilmanifolds.
 Let  $\liet=\lieg/\liec$.
 Below are some facts shown in Sections 2 and 3 of \cite{MPPS-2-step}.
 Since $\lieg$ is 2-step nilpotent, $\liet$ is abelian.
 As a vector space, $\lieg^{1,0}=\liet^{1,0}\oplus \liec^{1,0}$,
 and  $\lieg^{*(1,0)}=\liet^{*(1,0)}\oplus \liec^{*(1,0)}$.
 The only non-trivial Lie brackets in $\lieg^{1,0}\oplus \lieg^{0,1}$
 are of the form
 $ [\liet^{1,0},\liet^{0,1}]\subset \liec^{1,0}\oplus \liec^{0,1}.$ We assume that
 $\dim_\CC\liec^{1,0}=1$.

Explicitly, there exists a real basis $\{X_k,JX_k:1\leq k\leq n\}$ for $\liet$ and
$\{Z, JZ\}$ a real basis for $\liec$. The corresponding complex bases
for $\liet^{1,0}$ and $\liec^{1,0}$ are respectively composed of the following elements:
\begin{equation}
T_k=\frac{1}{2}(X_k-iJX_k) \quad \mbox{ and } \quad V=\frac{1}{2}(Z-iJZ).
\end{equation}
The structure equations of $\lieg$ are determined by
\begin{equation}\label{structure eq}
\lbra{\overline{T}_k}{T_j}=E_{kj} V-{\overline{E}}_{jk}{\overline{V}}
\end{equation}
for some constants $E_{kj}$.
Let $\{\omega^k: 1\leq k\leq n\}$ be the dual basis for $\liet^{*(1,0)}$, and let
$\{\rho\}$ be the dual basis for $\liec^{*(1,0)}$.
The dual structure equations for (\ref{structure eq}) are
\begin{equation}\label{dual 1}
d\rho=\sum_{i,j}E_{ji}\omega^i\wedge\oomega^j \quad \mbox{ and }
\quad d\omega^k=0.
\end{equation}
Equivalently,
\begin{equation}\label{dual 2}
d\overline{\rho}=-\sum_{i,j}\overline{E}_{ji}\omega^j\wedge\oomega^i
\quad \mbox{ and }
\quad d\oomega^k=0.
\end{equation}
It follows that
\begin{equation}\label{adj-T-on-rho-bar}
\lbra{T_j}{\overline{\rho}}={\cal L}_{T_j}{\overline{\rho}}
=\iota_{T_j}d{\overline{\rho}}=-\sum_i\overline{E}_{ji}\oomega^i.
\end{equation}

By Cartan formula (\ref{algebroid differential}),
\begin{equation}
\dbar {T_j}=\sum_{k}E_{kj}\oomega^k\wedge V
= \left(\sum_{k}E_{kj}\oomega^k\right)\wedge V.
\end{equation}

\

\noindent{\bf Example 1.} Consider a one-dimensional central extension of the Heisenberg
algebra $\lieh_{2n+1}$ of real dimension $2n+1$. Let $\{X_j, Y_j, Z, A: 1\leq j\leq n\}$
be basis with structure equations
\begin{equation}
\lbra{X_j}{Y_j}=-\lbra{Y_j}{X_j}=Z, \quad \mbox{ for all } \quad 1\leq j\leq n.
\end{equation}
The real center $\liec$ is spanned by $Z$ and $A$. Define an almost complex structure by
\[
JX_j=Y_j, \quad JY_j=-X_j, \quad JZ=A, \quad JA=-Z.
\]
It is an abelian complex structure. Let $V=\frac12(Z-iA)$ and
$T_j=\frac12(X_j-iY_j)$, then the complex structure
equations become
\[
\lbra{{\overline{T}_j}}{T_j}=-\frac{i}{2}(V+\overline{V}).
\]
Therefore, $E_{jj}=-\frac{i}{2}=-\overline{E}_{jj}$. Hence $d\rho$ is non-degenerate and serves
as an example for Theorem \ref{non degenerate theorem} and Proposition \ref{main proposition}.

\

\noindent{\bf Example 2.} On the direct sum of two Heisenberg
algebras $\lieh_{2m+1}\oplus \lieh_{2n+1}$, choose
a basis $\{X_j, Y_j, Z, A_k, B_k, C; 1\leq j\leq m, 1\leq k\leq n\}$
so that the non-zero structure equations are given by
\begin{equation}
\lbra{X_j}{Y_j}=-\lbra{Y_j}{X_j}=Z, \quad \lbra{A_k}{B_k}=-\lbra{B_k}{A_k}=C.
\end{equation}
Define an almost complex structure $J$ by
\[
JX_j=Y_j, \quad JA_k=B_k, \quad JZ=C.
\]
In fact, this defines an abelian complex structure. Let
\[
V=\frac12(Z-iC), \quad S_j=\frac12(X_j-iY_j), \quad T_k=\frac12(A_k-iB_k),
\]
so $\liec^{1,0}$ is spanned by $V$.
It follows that the non-zero complex structure equations are given as below
\[
\lbra{\overline{S}_j}{S_j}=-\frac{i}2(V+{\overline{V}}), \quad
\lbra{\overline{T}_k}{T_k}=\frac12(V-{\overline{V}}).
\]
We then have the structure constants
\[
E_{jj}=-\frac{i}2=-{\overline{E}_{jj}},
\quad \mbox{ and } \quad
E_{kk}=\frac12 ={\overline{E}_{kk}}
\]
for all $1\leq j\leq m$ and $1\leq k\leq n$.
It is now obvious that $d\rho$ is non-degenerate and serves as example for
Theorem \ref{non degenerate theorem} and Proposition \ref{main proposition}.

\

\noindent {\bf Example 3.} Consider a real vector space $P_{4n+2}$ spanned by
\[
\{X_{4k+1}, X_{4k+2}, X_{4k+3},
X_{4k+4}, Z_1, Z_2; 0\leq k\leq n-1\}.\] Define a Lie bracket by
\begin{eqnarray*}
\lbra{X_{4k+1}}{X_{4k+2}}=-\frac12Z_1,
\quad
\lbra{X_{4k+1}}{X_{4k+4}}=-\frac12Z_2,
\quad
\lbra{X_{4k+2}}{X_{4k+3}}=-\frac12Z_2.
\end{eqnarray*}
Define an abelian complex structure $J$ on  $P_{4n+2}$ by
\[
J{X_{4k+1}}={X_{4k+2}},
\quad
J{X_{4k+3}}=-{X_{4k+4}},
\quad
J{Z_1}=-Z_2,
\]
and define
\[
V=\frac12(Z_1+iZ_2), \quad
T_{2k+1}=\frac12 (X_{4k+1}-iX_{4k+2}),
\quad
T_{2k+2}=\frac12 (X_{4k+3}+iX_{4k+4}).
\]
The non-zero structure equations in terms of complex vectors are
\[
\lbra{{\overline{T}}_{2k+1}}{T_{2k+1}}=\frac{i}4(V+\overline{V});
\quad
\lbra{{\overline{T}}_{2k+1}}{T_{2k+2}}=-\frac14(V-\overline{V}).
\]
Hence
\[
E_{2k+1, 2k+1}=\frac{i}4,
\quad
E_{2k+1, 2k+2}=-\frac14,
\quad
E_{2k+2, 2k+1}=-\frac14.
\]
It follows that $d\rho$ is non-degenerate, and hence (\ref{dbar exact equation}) is solvable,
providing another example for Theorem \ref{non degenerate theorem} and Proposition \ref{main proposition}.

\

\noindent{\bf Example 4.} Consider a real vector space $W_{4n+6}$ spanned by \[
\{X_{4k+1}, X_{4k+2}, X_{4k+3},
X_{4k+4}, Z_1, Z_2; 0\leq k\leq n\}.\] Define a Lie bracket by
\begin{eqnarray*}
\lbra{X_{4k+1}}{X_{4k+3}}=-\frac12Z_1,
& \quad &
\lbra{X_{4k+1}}{X_{4k+4}}=-\frac12Z_2, \\
\lbra{X_{4k+2}}{X_{4k+3}}=-\frac12Z_2,
& \quad &
\lbra{X_{4k+2}}{X_{4k+3}}=\frac12Z_1.
\end{eqnarray*}
One can define an abelian complex structure $J$ by
\[
JX_{4k+1}=X_{4k+2}, \quad JX_{4k+3}=-X_{4k+4}, \quad JZ_1=-Z_2.
\]
For $0\leq k\leq n$, define
\[
V=\frac12(Z_1+iZ_2), \quad
T_{2k+1}=\frac12 (X_{4k+1}-iX_{4k+2}),
\quad
T_{2k+2}=\frac12 (X_{4k+3}+iX_{4k+4}).
\]
It is now a straightforward computation to show that the non-zero structure equations are
\begin{equation}
\lbra{{\overline{T}}_{2k+1}}{T_{2k+2}}=-\frac12 V, \quad
\lbra{{\overline{T}}_{2k+2}}{T_{2k+1}}-=\frac12{\overline{V}}.
\end{equation}
Except when
\[
E_{2k+1,2k+2}=-\frac12, \quad \mbox{ for all } \quad 0\leq k\leq n,
\]
all other structure constants are equal to zero. In particular,
\begin{equation}
d\rho=-\frac12\sum_{k=0}^n\omega^{2k+2}\wedge\oomega^{2k+1}, \quad
\mbox{ and }
\quad
d\overline{\rho}=\frac12 \sum_{k=0}^n\omega^{2k+1}\wedge\oomega^{2k+2}.
\end{equation}
Moveover,
\begin{equation}\label{holomorphic vector fields}
\dbar T_{2k+1}=0, \quad \dbar T_{2k+2}=-\frac12 \oomega^{2k+1}\wedge V.
\end{equation}

Treating $d\rho$ and $d\overline{\rho}$ as maps from $\liet^{1,0}$ to
$\liet^{0,1}$, their image spaces are transversal.
Therefore, given $T\in \liet^{1,0}$ such that there exists $X\in \liet^{1,0}$ with
$\iota_Td\overline{\rho}=-\iota_Xd{\rho}$ only if $\iota_Td\overline{\rho}=0$. It is possible only when
$T$ is the complex linear span of
$\{ T_2, \dots, T_{2n+2} \}$. In this case,
it simply means that $\lbra{T}{\overline\rho}=0$, and hence $\adL{\overline\rho}=0$.
This example illustrates the
situation in Proposition \ref{secondary proposition} as well as
the non-trivial conditions in Theorem \ref{non degenerate theorem}.

It is also apparent that if we choose $\Lambda=V\wedge T_{2k+1}$, then
\[
\adL(\overline\rho)
=V\wedge \iota_{T_{2k+1}}d\overline{\rho}=-\frac12 V\wedge \oomega^{2k+2},
\]
which is not $\dbar$-exact. Therefore, the map $\phi$ in Theorem \ref{phi inj}
fails to be surjective and the holomorphic
Poisson spectral sequence fails to degenerate on
its first page.

In this example, if we vary $T$ from $\{ T_2, \dots, T_{2n+2} \}$
to $\{ T_1, \dots, T_{2n+1} \}$ and define $\Lambda=V\wedge T$,
the dimension of the cohomology
$H^{\bullet}_{\Lambda}$ jumps down.

Consider $\Lambda=V\wedge T$ where $T$ is in the linear span of
$\{ T_2, \dots, T_{2n+2} \}$. The map $\phi$ is then an identity map.
In particular,
\begin{equation}
H^1_{\Lambda}=H^0(\lieg^{1,0})\oplus H^1(\lieg^{0,0}).
\end{equation}
From the structural equations (\ref{holomorphic vector fields}),
\begin{equation}
H^0(\lieg^{1,0})= \{V, T_1, \dots, T_{2n+1} \}, \quad
H^1(\lieg^{0,0})= \lieg^{*(0,1)}.
\end{equation}
Therefore,
\begin{equation}\label{dim = 3n+5}
\dim H^1_{\Lambda}=\dim H^0(\lieg^{1,0}) +\dim H^1(\lieg^{0,0})=(n+2)+(2n+3)=3n+5.
\end{equation}

\noindent {\bf A Further Observation.} We can consider the virtual parameter space of generalized complex deformations, utilizing the previous example.  For the second hypercohomology, we have the Hodge-like decomposition
\begin{equation}
H^2_{\Lambda}=H^0(\lieg^{2,0})\oplus H^1(\lieg^{1,0})\oplus H^2(\lieg^{0,0}).
\end{equation}

We could consider $H^0(\lieg^{2,0})=H^0(M, \Theta^{2,0})$ as deformations of
holomorphic Poisson structures by varying $\Lambda$ in the space of holomorphic
bivector fields without changing the complex structure. Moreover, we could
consider $H^1(\lieg^{1,0})=H^1(M, \Theta^{1,0})$ to be the virtual parameter space
of classical deformation of complex structure. It remains to examine
the nature of $H^2(\lieg^{0,0})=H^2(M, {\calo})$.

By Lemma \ref{0q forms}
\begin{equation}
H^2(\lieg^{0,0})=\lieg^{*(0,2)}
=\liec^{*(0,1)}\otimes \liet^{*(0,1)}\,\oplus\, \liet^{*{0,2}}.
\end{equation}

Note that every element in $\liet^{*{0,2}}$ is not only $\dbar$-closed but also $d$-closed.
Therefore, the deformations that they generate are  $B$-field transformations of
$\Lambda$ as a generalized complex structure. In terms of generalized complex geometry, they are equivalent
to the given holomorphic Poisson structure $\Lambda=V\wedge T$
\cite{Marco}.
However, elements in $\liec^{*(0,1)}\otimes \liet^{*(0,1)}$ are $\dbar$-closed but not $d$-closed.
For instance, consider $\Oomega=\overline{\rho}\wedge \oomega^{2h+1}$ where $\omega^{2h+1}$
is dual to $T_{2h+1}$. As it satisfies both $\dbar_{\Lambda}\Oomega=0$ and
$\lbra{\Oomega}{\Oomega}=0$, it defines an integrable generalized deformation of $\Lambda$
\cite{Marco}. The corresponding
Lie algebroid differential after deformation is  given by
\begin{equation}
\overline{\delta}=\dbar_\Lambda+\lbra{\Oomega}{-}.
\end{equation}
As both $\Lambda=V\wedge T$ and $\oomega^{2h+1}$ are in the center of Schouten algebra $\oplus_{p,q}B^{p,q}$,
\begin{equation}
\overline{\delta}=\dbar-\oomega^{2h+1}\wedge\lbra{\overline\rho}{-}.
\end{equation}
Therefore, every element in $\lieg^{*(0,1)}$ is in the kernel of $\overline\delta$.
Next we examine the action of $\overline\delta$ on $\lieg^{1,0}$.
Apparently, ${\overline\delta}V=0$. For each $k$, according to (\ref{holomorphic vector fields})
\begin{eqnarray*}
\overline{\delta}T_{2k+1} &=&\dbar T_{2k+1}-\oomega^1\wedge\lbra{\overline\rho}{T_{2k+1}}
=\oomega^{2h+1}\wedge\lbra{T_{2k+1}}{\overline\rho}=\frac12\oomega^{2h+1}\wedge \oomega^{2k+2}\\
\overline{\delta}T_{2k+2} &=&\dbar T_{2k+2}-\oomega^{2h+1}\wedge\lbra{\overline\rho}{T_{2k+2}}
=\dbar T_{2k+2}=-\frac12 \oomega^{2k+1}\wedge V.
\end{eqnarray*}
In particular, on the first cohomology level
\begin{equation}
\ker \overline{\delta}=\liec^{1,0}\oplus \lieg^{*(0,1)}.
\end{equation}
In view of the result in (\ref{dim = 3n+5}),
the first cohomology changes after deformation. The class $\Oomega=\overline\rho\wedge \oomega^{2h+1}$
defines a non-trivial generalized complex deformation of $\Lambda=V\wedge T$.

\

\noindent{\bf Acknowledgment.}
Y.~S. Poon thanks the Institute of Mathematical Sciences at the Chinese University of Hong Kong for
hosting his numerous visits in the past two years.

\noindent{Yat Sun  Poon:
    Department of Mathematics, University of California,
    Riverside, CA 92521, U.S.A.. E-mail: ypoon@ucr.edu.}

\

\noindent{John Simanyi:
    Department of Mathematics, University of California,
    Riverside, CA 92521, U.S.A.. E-mail: jsima003@ucr.edu.}

\end{document}